\documentclass[11pt,reqno]{amsproc}
\usepackage[german,english]{babel}  
\usepackage[latin1]{inputenc} 
\usepackage{graphicx, color}
\usepackage{amssymb}
\numberwithin{equation}{section}
\newtheorem{theorem}{Theorem}[section]

\newtheorem{lemma}[theorem]{Lemma}

\newtheorem{remark}[theorem]{Remark}
\newtheorem{remarks}[theorem]{Remarks}
\newcommand{\Z}{\mathbb{Z}}

\renewcommand{\P}{\mathbb{P}}
\newcommand{\E}{\mathbb{E}}
\newcommand{\e}{\mathcal{E}}

\newcommand{\eh}{\frac{1}{2}}

\newcommand{\gf}[1]{`#1'}

\def\<{\langle}
\def\>{\rangle}

\def\p0{\psi_0}                              

\def\1Ll{{1 \over {\vert \Lambda_L \vert}}}

\def\j1d{{(1+\vert j \vert )}^{-(d+2)}}
\def\x1a{{(1+ \vert x\vert )}^{-\alpha}}
\def\2j1a{{(1+ \vert j\vert )}^{-\alpha}}

\begin{document}

\title[square-root law]
{On Penrose's square-root law and beyond}

\author{Werner Kirsch}
\address{Institut f\"ur Mathematik,
Ruhr-Universit\"at Bochum,\newline\indent D-44780 Bochum, Germany\vspace{2mm}}
\email{werner.kirsch@ruhr-uni-bochum.de}

\maketitle

\begin{abstract}
In certain bodies, like the Council of the EU, the member states
have a voting weight which depends on the population of the
respective state. In this article we ask the question which voting
weight guarantees a `fair' representation of the citizens in the
union. The traditional answer, the square-root law by Penrose, is
that the weight of a state (more precisely: the voting power) should
be proportional to the square-root of the population of this state.
The square root law is based on the assumption that the voters in
every state cast their vote independently of each other. In this
paper we concentrate on cases where the independence assumption is
not valid.
\end{abstract}

\section{Introduction\label{sec:intro}}
All modern democracies rely on the idea of representation.
A certain body of representatives, a parliament for example, makes
decisions on behalf of the voters. In most parliaments each of its members represents
roughly the same number of people, namely the voters in his or her
constituency.

There are other bodies in which the members represent different
numbers of voters. A prominent example is the Council of the
European Union. Here ministers of the member states represent the
population of their respective country. The number of people
represented in the different states  differs from about 400,000 for
Malta to more than 82 million for Germany. Due to this fact the
members of the Council have a certain number of votes depending on
the size of the country they represent, e.g. 3 votes for Malta, 29
votes for Germany. The votes of a country cannot be split, but have
to be cast as a block.\footnote{The current voting system in the
Council is based on the treaty of Nice. It has additional components
to the procedure described above, which are irrelevant in the
present context. For a description of this voting system and further
references see e.g. \cite{Kirsch}.}

Similar voting systems occur in various other systems, for example
in the Bundesrat, Germany's state chamber of parliament and in the
electoral college in the USA.\footnote{The electoral college is not
exactly a heterogeneous voting system in the sense defined below,
but it is very close to it.}

Let us call such a system in which the members represent subsystems
(states) of different size a \emph{heterogeneous voting system}. In the
following we will call the assembly of representatives in a
heterogeneous voting system the \emph{council}, the sets of voters
represented by the council members the \emph{states}.

It is quite clear, that in a heterogeneous voting system a bigger
state (by population) should have at least as many votes in the
council as a smaller state. It may already be debatable whether the
bigger states should have \emph{strictly} more votes than the smaller
states (cf. the Senate in the US constitution). And if yes, how much more votes the bigger state should get?

In this note we address the question: `What is a fair distribution
of power in a heterogeneous voting system?'

There exist various answers to this question, depending on the
interpretation of the words `fair' and `power'.

The usual and quite reasonable way to formulate the question in an
exact way is to use the concept of power indices. One calls a
heterogeneous voting system fair if all voters in the member states
have the same influence on decisions of the council. By `same
influence' we mean that the power index of each voter is the same
regardless of her or his home state. If we choose then Banzhaf power
index to measure the influence of a voter we obtain the celebrated
Penrose's square-root law (see e.g. \cite{FMbook}).

The square-root law states that the distribution of power in a
heterogeneous voting system is fair if the power (index) of each
council member $i$ is proportional to $\sqrt{N_i}$, where $N_i$ is
the population of the state which $i$ represents.

In their book \cite{FMbook} Felsenthal and Machover formulate a
second square-root law. There they base the notion of
`fairness' on the concept of \emph{majority deficit}.

The majority deficit is zero if the voters favoring the decision of
the council are the \emph{majority}. If the voters favoring the
decision of the council are the \emph{minority} then the majority
deficit is the margin between the number of voters objecting to the
decision and those agreeing with it (see Def. 3.3.16 in
\cite{FMbook}).

The notion of fairness we propose in this paper is closely related to the
concept of majority deficit. We will call a decision of the council
\emph{in agreement with the popular vote} if the percentage of
voters agreeing with a proposal (popular vote) is as close as
possible to the percentage of council votes in favor of the
proposal. (We will make this notion precise in the next section.)

For both concepts we have to average over the possible voting
configurations. This is usually done by assuming that voters vote
independently of each other. The main purpose of this note is to
investigate some (we believe reasonable) models where voters do not
vote independently.

We will discuss two voting models with voting behavior which
is \emph{not} independent. The first model considers societies which
have some kind of \gf{common belief}. A typical situation of this kind
is a strong religious group (or church) influencing the voting
behavior of the voters. This model is discussed in detail in Section  \ref{sec:cbm}.

In the other model voters tend to vote the same way \gf{the majority
does}. This is a situation where voters do not want to be different
from others. We call this the \emph{mean field model} referring to
an analogous model from statistical physics. See Section \ref{sec:mfm} for
this model.

In fact, both models can be interpreted in terms of statistical
physics. Statistical physics considers (among many other things)
magnetic systems. The elementary magnet, called a spin, has two
possible states which are \gf{$+1$} or \gf{$-1$} (spin up, spin down).
This models voting \gf{yes} or \gf{no} in a voting system. Physicists
consider different kinds of interactions between the single spins,
one given through an exterior magnetic-field - corresponding to a
society with \gf{a common belief} - or through the tendency of the
spins to align - corresponding to the second voting model. We discuss
the analogy of voting models with spin systems in Section \ref{sec:spin}.

Our investigations of voting models with statistical dependence is
much inspired by the paper \cite{LarVal}. The first model is also
based on the work by Straffin \cite{Straffin}.

It does not come as a surprise that we obtain a square-root law for
a model with independent voters, just as in the case considered by
Felsenthal and Machover (\cite{FMbook}).

For the mean field model we still get a square-root law for the best
possible representation in the council \emph{as long as} the mutual
interaction between voters is not too strong.

However as the coupling between voters exceeds a certain threshold,
the fairest representation in the council is no longer given by
votes proportional to $\sqrt{N_i}$ but rather by votes proportional to
$N_i$. This is a typical example of a phase transition.

In the model of common belief the fair representation weight depends on
the strength of the common belief for large populations. If this strength
is independent of the population size fair representation is almost
always given by voting weights proportional to $N_i$, the square-root law occurring only in marginal cases.
However, if the common belief
decreases with increasing population one can get any power law behavior ${N_i}^\alpha$
for the optimal weight as long as $\eh\leq \alpha\leq 1$. In fact, statistical investigations
on real life data suggest that this might happen (see \cite{GelKB}).

We leave the mathematical proofs of our results for the appendices (Sections \ref{sec:pr} to \ref{sec:pr2}).

\medskip
\textbf{Acknowledgment:} It is a pleasure to thank Hans-J\"urgen Sommers, Duisburg-Essen and Wojciech S\l omczy\'nski
and Karol \.Zyczkowski, Krakow for valuable discussions.

\section{The general model \label{sec:model}}

We consider $N$ voters, denoted by $1, 2, \ldots, N$. Each of them
may vote `yes' or `no'; abstentions are not allowed. The vote of the
voter $i$ is denoted by $X_{i}$.

The possible voting results are $X_{i}=+1$ representing \gf{yes}
and $X_{i}=-1$ for \gf{no}. We consider the quantity $X_{i}$ as
random, more precisely there is a probability measure $\mathbb{P}$
on the space $\{ -1,1 \}^N$ of possible voting results. This measure
will be specified later. The conventional assumption on $\mathbb{P}$
is that the random quantities $X_{i}$ are independent from each
other, but we are \emph{not} making this assumption here.

Our interpretation of this model is as follows. The voters react on
a proposal in a rational way, that is to say: A voter does
\emph{not} roll a dice to determine his or her voting behavior but
he or she votes for or against a given proposal according to his/her
personal belief, knowledge, experience etc. It is rather the
proposal which is the source of randomness in this system. We
imagine the voting system is fed with propositions in a completely
random way. This could be either a real source of proposals or just
a Gedankenexperiment to measure the behavior of the voting system.

The rationality of the voters implies that a voter who casts a
\gf{yes} on a certain proposition will necessarily vote \gf{no} on
the diametrically opposed proposition. Since we assume that the
proposals are completely random any proposal and its antithetic
proposal must have the same probability. This implies
\begin{eqnarray}\label{eq:einhalb}
\mathbb{P}(X_{i}=1)=\mathbb{P}(X_{i} = -1) = \frac{1}{2} \ .
\end{eqnarray}

More generally, we conclude that
\begin{eqnarray}\label{eq:symm}
\mathbb{P}(X_{i_1}=\xi_1 , ... ,X_{i_r}=\xi_r ) =
\mathbb{P}(X_{i_1}=-\xi_1 , ... ,X_{i_r}=-\xi_r )
\end{eqnarray}
for any set $i_1...,i_r$ of voters and any
$\xi_1,...\xi_r \in \{ -1,1 \}$.

We call the property (\ref{eq:symm}) the \emph{symmetry} of the voting
system. Any measure $\P$ satisfying (\ref{eq:symm}) is called a
\emph{voting measure}.

The symmetry assumption (\ref{eq:symm}) does \emph{not} fix the probability measure
$\mathbb{P}$. Only if we assume in addition that the $X_{i}$ are
statistically independent we can conclude from (\ref{eq:symm}) that

\begin{eqnarray}\label{eq:symmind}
\mathbb{P}(X_{i_1}=\xi_1 , ... ,X_{i_r}=\xi_r ) =
(\frac{1}{2})^r \ .
\end{eqnarray}


So far, we have not specified any decision rule for the voting
system. The above probabilistic setup is completely independent from
the voting rule, a fact which was emphasized in the work
\cite{LarVal}.

A \emph{simple majority rule} for $X_1,\dots,X_N$ is given by the
decision rule: Accept a proposal if $\sum_{j=1}^N X_j>0$ and reject
it otherwise.

By a \emph{qualified majority rule} we mean that at least a
percentage $q$ (called the \emph{quota}) of votes is required for
the acceptance of a proposal. In term of the $X_j$ this means:

\begin{equation}\label{def:qm}
\sum_{j=1}^N X_j\geq(2q-1)N.
\end{equation}

Indeed, it is not hard to see that the number of affirmative votes
is given by

\begin{equation}\nonumber
\frac{1}{2}\left(\sum_{j=1}^N X_j+N\right).
\end{equation}

From this the assertion (\ref{def:qm}) follows.

In particular, the simple majority rule is obtained form
(\ref{def:qm}) by choosing $q$ slightly bigger than $\frac{1}{2}$.

The sum $\sum_{j=1}^N X_j$ gives the difference between the number of
\gf{yes}-votes and the number of \gf{no}-votes. We call the quantity

\begin{equation}\label{def:margin}
M(X):=\left|\sum_{j=1}^N X_j\right|
\end{equation}

the margin of the voting outcome $X=(X_1,\dots,X_N)$. It measures
the size of the majority with which the proposal is either accepted
or rejected in simple majority voting.

In qualified majority voting with quota $q$ the corresponding
quantity is the \emph{$q$-margin} $M_q(X)$ given by:

\begin{equation}\label{def:qmargin}
M_q(X):=\left|\sum_{j=1}^N X_j-(2q-1)N\right|.
\end{equation}

Now, we turn to voting in the council. We consider $M$ states, the state number $\nu$ having $N_\nu$ voters.
Consequently the total number of voters is $N=\sum N_\nu$. The
vote of the voter $i$ in state $\nu$ is denoted by $X_{\nu i}$,
$\nu = 1,...,M$ and $i=1,...,N_\nu$.\footnote{We label the states
using Greek characters and the voters within a state by Roman
characters.}

We suppose that each state
government knows the opinion of (the majority of) the voters in that state and
acts accordingly.\footnote{Although this is the central idea of
representative democracy this idealization may be a little naive in
practice.} That is to say: If the majority of people in state $\nu$
supports a proposal, i.e. if

\begin{equation}\label{affstvote}
\sum_{i=1}^{N_\nu}X_{\nu i}>0
\end{equation}

then the representative of state $\nu$ will vote \gf{yes} in the
council otherwise he or she will vote \gf{no}. If we set $\chi(x)=1$
for $x>0$, $\chi(x)=-1$ for $x\leq0$ the representative of state
$\nu$ will vote

\begin{equation}\label{def:xi}
\xi_\nu=\chi\left(\sum_{i=1}^{N_\nu}X_{\nu i}\right)
\end{equation}

in the council. If the state $\nu$ has got a weight $w_\nu$ in the
council the result of voting in the council is given by:

\begin{equation}\label{votecounc}
\sum_{\nu=1}^M
w_\nu\,\xi_\nu=\sum_{\nu=1}^Mw_\nu\;\chi\!\left(\sum_{i=1}^{N_\nu}X_{\nu
i}\right).
\end{equation}

Thus, the council's decision is affirmative if
$\sum_{\nu=1}^Mw_\nu\xi_\nu$ is positive, provided the council votes
according to simple majority rule.

The result of a popular vote in all countries $\nu=1,\dots,N$ is

\begin{equation}\label{popvote}
P~=~\sum_{\nu=1}^M\sum_{i=1}^{N_\nu}X_{\nu i}.
\end{equation}

We will call voting weights $w_\nu$ for the council \emph{fair} or \emph{optimal}, if the
council's vote is as close as possible to the public votes. To make
this precise let us define

\begin{equation}
C=\sum_{\nu=1}^Mw_\nu\,\chi\!\left(\sum_{i=1}^{N_\nu}X_{\nu
i}\right)
\end{equation}

\noindent the result of the voting in the council. Both $P$ and $C$ are random
quantities which depend on the random variables $X_{\nu i}$. So, we
may consider the mean square distance $\vartriangle$ between $P$ and
$C$, i.e. denoting the expectation over the random quantities by
$\E$, we have

\begin{eqnarray}
\vartriangle &=& \E\left((P-C)^2\right)\\ \label{msdelta}&=&
\E\left(\left\{\sum_{\nu=1}^M\sum_{i=1}^{N_\nu}X_{\nu i}-\sum_{\nu=1}^M
w_\nu\chi\big(\sum_{i=1}^{N_\nu}X_{\nu i}\big)\right\}^2\right).
\end{eqnarray}

In a democratic system the decision of the council should be as
close as possible to the popular vote, hence we call a system of
weights \emph{fair} or \emph{optimal} if
$\vartriangle=\vartriangle(w_1,\dots,w_M)$ is minimal among all
possible values of $w_\nu$.

In the following we suppose that the random variables $X_{\nu i}$
and $X_{\mu j}$ are independent for $\nu\neq\mu$. This means that
voters in different states are not correlated. We do not assume at
the moment that two voters from the same state vote independently of each
other.

We have the following result:

\begin{theorem}\label{th:srl}
Fair voting in the council is obtained for the values

\begin{eqnarray*}
w_\nu&=&\,\eh~\E\left(\left|\sum_{i=1}^{N_\nu}X_{\nu
i}\right|\right)\\
&=&\,\eh~\E\Big(M(X_{\nu 1},\dots,X_{\nu
{N_{\nu}}})\Big).
\end{eqnarray*}
\end{theorem}

This result can be viewed as an extension of Penrose's square-root
law to the situation of correlated voters. We will see below that it
gives $w_\nu\sim\sqrt{N_\nu}$ for independent voters.

Theorem \ref{th:srl} has a very easy - we hope convincing -
interpretation: $w_\nu$ is the expected margin of the voting result
in state $\nu$. In other words, it gives the expected number of
people in state $\nu$ that agree with the voting of $\nu$ in their
council minus those that disagree, i.e. the net number of voters
which the council member of $\nu$ actually represents.

If we choose any multiple $cw_1,\dots,cw_{N_\nu}$ $(c>0)$ of the
weights $w_1,\dots.w_{N_\nu}$ we obtain the same voting system as
the one defined by $w_1,\dots,w_n$. In this sense the weight $w_\nu$
of Theorem \ref{th:srl} are not unique, but the voting system is.

We will prove Theorem \ref{th:srl} in section \ref{sec:pr}. We
remark that the proof requires the symmetry assumption
(\ref{eq:symm}) and the independence of voters from \emph{different}
states.

The next step is to compute the expected margin $\E(M(X_{\nu}))$, at
least asymptotically for large number of voters $N_\nu$. This quantity
depends on the correlation structure between the voters in state
$\nu$. As we will see, different correlations between voters give
very different results for $\E(M(X_{\nu}))$ and hence for the
optimal weight $w_\nu$.

We begin with the classical case of independent voters.

\begin{theorem}\label{th:indep}
If the voters in state $\nu$ cast their votes independently of each
other then

\begin{equation}
\E\left(\left|\sum_{i=1}^{N_\nu}X_{\nu i}\right|\right)\sim
c\,\sqrt{N_\nu}
\end{equation}

for large $N_\nu$.
\end{theorem}

Thus, we recover the square-root law as we expected. (For the square-root law see Felsenthal and Machover \cite{FMbook}.) In terms of
power indices the independence assumption is associated to the
Banzhaf power index. Therefore, it is not surprising that also the
Banzhaf index leads to a square-root rule.

It is questionable (as we know from the work of Gelman, Katz and Bafumi \cite{GelKB})
whether the independent voters model is valid in many real-life voting systems.
This is one of the reasons to extend the model as we do in the present paper.

\section{The `common belief' model\label{sec:cbm}}

In this section we consider a model we dub the `common belief
model'. It generalizes a voting measure introduced and investigated
by Straffin \cite{Straffin} in connection with the Shapley-Shubik
power index.

We imagine that inside a certain society there is a strong common
belief which may for example be due to a powerful religious group,
a generally accepted political ideology or a strong tradition. This causes a tendency to a creation of
strong majorities in a certain type of questions. For example, in a
country with a strong catholic majority there may be a strongly
correlated view about abortion among voters, but, may be, not about
speed limits on highways. One might have a similar effect if a
person dominates the public or private media or both.

We model such a situation by introducing a random variable $Z$ which
reflects the `common belief' on the subject at hand. $Z=+1$ means
that all voters agree to accept the given proposal, $Z=-1$ means
that all voters will reject it. The random variable $Z$ is allowed
to take any value in $[-1 , 1]$. If $Z=0$ there is no common belief
on the proposal. If $Z>0$ there is some common belief favoring the
proposal which is weak if $Z$ is close to 0 and strong if it is
close to 1. The probability distribution of $Z$ is denoted by $\mu$,
hence
\begin{equation}\label{def:distZ}
\mu ([a,b])=\P(Z\in [a,b]).
\end{equation}

$Z$ has to satisfy a symmetry
condition similar to (\ref{eq:einhalb}), namely
\begin{equation}
\mathbb{P}(Z\in[a, b]) = \mathbb{P}(Z\in[-b, -a]),
\end{equation}
i.e.
\begin{equation}\label{eq:symmZ}
\mu([a, b]) = \mu([-b, -a]) \ .
\end{equation}

In our model the \gf{common belief} variable $Z$ influences the
probability that a voter $i$ votes $\pm 1$. Given $Z=\zeta \in
[-1,1]$, then the conditional probability that $X_i=1$ given $Z=\zeta$
is given by:
\begin{eqnarray}\label{eq:cbm1}
\mathbb{P}(X_i = 1 | Z=\zeta) = \frac{1}{2}(\zeta + 1)~=~p_\zeta \ .
\end{eqnarray}
Thus, if $\zeta = 1$ any voter $i$ will vote "$+1$" with probability
1, if $\zeta = -1$ he or she will vote  "$+1$" with probability 0
and if $\zeta = 0$ then $i$ votes with probability one half in favor
or against the proposal.

We denote by $P_s$ the probability measure on $\{-1, +1\}^N$ with
$P_s(X_i=1)=s$ and $P_s(X_i=-1)=1-s$, which makes the $X_i$ independent.
We also denote by $E_s$ the expectation with respect to $P_s$.
 Note that the probability $p_\zeta=\frac{1}{2}(1+\zeta)$ in (\ref{eq:cbm1}) is
chosen in such a way that $E_{p_\zeta}(X_i)=\zeta$. Thus the value
of the \gf{common belief} variable $Z$ gives the expected voting
result of a single voter.

For \emph{given} $Z=\zeta$ we assume the $X_i$ to be independent. Thus we
have
\begin{eqnarray}\label{eq:cbmfull}
\mathbb{P}(X_1 = \xi_1 , ... , X_N = \xi_N) = \int ~\Big(\;\prod_1^N
\,P_{p_\zeta}(X_i = \xi_i)\,\Big)\; d\mu (\zeta) \ .
\end{eqnarray}

The measure
$\mathbb{P}$ in (\ref{eq:cbmfull}) depends on the probability
distribution $\mu$, hence we sometimes denote it by
$\mathbb{P}_\mu$.

The probability measure $\mathbb{P}_\mu$ defined in
(\ref{eq:cbmfull}) satisfies the symmetry condition (\ref{eq:symm})
due to assumption (\ref{eq:symmZ}). Of course, $\mathbb{P}_\mu$
defines a whole class of examples, each (symmetric) probability
measure $\mu$ on $[-1,1]$ defines its unique $\mathbb{P}_\mu$. If we
choose $\mu = \delta_0$, i.e. $\mu([a,b])=1$ if $a\leq 0 \leq b$ and
$=0$ otherwise, we obtain independent random variables $X_i$ as
discussed in the final part of section \ref{sec:model}. Indeed,
$\mu = \delta_0$ means that $Z=0$, consequently (\ref{eq:cbmfull})
defines independent random variables. Observe, that this is the only
measure for which $Z$ assumes a fixed value ($\mu$ has to be
symmetric!).

Another interesting example is the case when $\mu$ is the uniform
distribution on $[-1,1]$. This case was considered by Straffin
\cite{Straffin}. He observed that this model is intimately connected
with the Shapley-Shubik power index.

To apply the \gf{common belief} model to a given heterogeneous
voting model we have to specify the measure $\mu$, of course. In
fact, this measure may change from state to state. In particular,
one may argue that larger states tend to have a less homogeneous
population and hence a weaker influence of a religious or political
group. For example, we will later discuss a model modifying
Straffin's example where $\mu (dz)=\frac{1}{2}\chi_{[-1,1]}(z)dz$ to
a measure where $\mu_N$ depends on the population $N$, namely
\begin{equation}\label{ex:Strafmod}
\mu_N (dz) = \frac{1}{2a_N}\chi_{[-a_N, a_N]}(z)dz
\end{equation}
with parameters $0<a_N\leq 1$. In particular, if we have $a_N \to 0$ as $N\rightarrow \infty$, the
parameter $a_N$ reflects the tendency of a common belief to decrease
with a growing population.

Except for the trivial case $\mu = \delta_0$ the random variables
$X_i$ are never independent under $\mathbb{P}_\mu$. This can be seen
from the covariance
\begin{eqnarray}\label{def:cov}
\langle X_i , X_j \rangle_\mu := \mathbb{E}_\mu(X_i
X_j)-\mathbb{E}_\mu(X_i)\mathbb{E}_\mu(X_j) \ .
\end{eqnarray}

In (\ref{def:cov}) as well as in the following $\mathbb{E}_\mu$
denotes expectation with respect to $\mathbb{P}_\mu$. In fact, the
random variables $X_i$ are always positively correlated:

\begin{theorem}\label{th:corrpmu}
For $i \neq j$ we have
\begin{eqnarray}\label{eq:corrpmu}
\langle X_i , X_j \rangle_\mu = \int \zeta^2 \,d\mu (\zeta) \ .
\end{eqnarray}
\end{theorem}

The quantity $\int \zeta^2 \,d\mu (\zeta)$ is called the second
moment of the measure $\mu$. Since the first moment $\int \zeta
\,d\mu (\zeta)$ vanishes due to (\ref{eq:symmZ}) the second moment
equals the variance of $\mu$. Observe that $\int \zeta^2 d\mu
(\zeta) = 0$ implies $\mu = \delta_0$. For independent random
variables $\langle X_i , X_j \rangle_\mu = 0$, so (\ref{eq:corrpmu})
implies that $X_i, X_j$ depend on each other unless $\mu =
\delta_0$.

To investigate the impact of the common belief measure $\mu$ on the
ideal weight in a heterogeneous voting model we have to compute the
quantity
\begin{eqnarray}\label{eq:weightmu}
\mathbb{E}_\mu (|\sum X_i|)
\end{eqnarray}
for a measure $\mu$ and population $N$ (at least for large $N$).
This is done with the help of the following Theorem:

\begin{theorem}\label{th:weightmu}
$\Big|\; \mathbb{E}_\mu(\frac{1}{N}\,| \sum_1^N X_i \;|) - \int\,
|\zeta|\;d\mu (\zeta) \Big| \leq \frac{1}{\sqrt{N}}$.
\end{theorem}

If we choose $\mu \neq \delta_0$ independent of the (population of
the) state Theorem \ref{th:weightmu} implies that the optimal weight
in the council is \emph{proportional} to $N$ (rather than
$\sqrt{N}$). This is true in particular for the original
Straffin model \cite{Straffin} where $\mu_n\equiv \eh \chi_{[-1,1]}(z)\,dz$
which corresponds to the Shapley-Shubik power index.

Let us define $\overline{\mu} = \int |\zeta|\,d\mu (\zeta)$. If $\mu =
\mu_N$ depends on the population then
\begin{displaymath}
\mathbb{E}_{\mu_N} ( | \sum X_i | ) \thicksim ~N\,\overline{\mu}_N
\end{displaymath}
as long as $\overline{\mu}_N \geq \frac{1}{N^{1/2-\varepsilon}}$ for
some $\varepsilon > 0$. However, if $\overline{\mu}_N \leq
\frac{1}{N^{1/2-\varepsilon}}$, then
\begin{displaymath}
E_{\mu_N}(| \sum X_i |) \thicksim \sqrt{N} \ .
\end{displaymath}
Hence, in this case we rediscover a square-root law.

We summarize:

\begin{theorem}\label{th:summu}
Let us suppose that a state with a population of size $N$ is
characterized by a common belief measure $\mu_N$, then:
\begin{enumerate}
\item If
    \begin{equation}\label{ass:geqsqrt}
    \overline{\mu}_N=\int |\zeta|\;d\mu_N(\zeta)\geq \,C\,\frac{1}{N^{1/2-\varepsilon}}
    \end{equation}
    for some $\varepsilon>0$ and for all large $N$ then  the optimal weight $w_N$ is given by:
    \begin{equation}\label{res:geqsqrt}
    w_N\;=\;\mathbb{E}_\mu(\,|\sum_1^N X_i
    |\,)\;\sim\;N\,\overline{\mu}_N.
    \end{equation}

\item If
    \begin{equation}\label{ass:leqsqrt}
    \overline{\mu}_N=\int |\zeta|\;d\mu_N(\zeta)\leq \,C\,\frac{1}{N^{1/2+\varepsilon}}
    \end{equation}
    then for large $N$ the optimal weight $w_N$ is given by:
    \begin{equation}\label{res:leqsqrt}
    w_N\;=\;\mathbb{E}_\mu(\,|\sum_1^N X_i |\,)\;\sim\;\sqrt{N}.
    \end{equation}

\end{enumerate}
\end{theorem}

\noindent\textbf{Example:} In our Straffin-type example
(\ref{ex:Strafmod}) we choose:
\begin{equation}
\mu_N (dz) = \frac{1}{2a_N}\chi_{[-a_n, a_n]}(z)dz,
\end{equation}
then:
\begin{equation}
\overline{\mu}_N  = \eh\,a_N.
\end{equation}

So, if $a_N\leq C\,\frac{1}{\sqrt{N}}$ we have $w_N\sim\sqrt{N}$, otherwise we obtain $w_N\sim a_N$.

\begin{remarks}\noindent\\[-5mm]
\begin{enumerate}
\item Our result shows that in all cases the optimal weight $w_N$
satisfies $C\,\sqrt{N}\leq w_N\leq N$. It is a matter of empirical
studies to determine which measure $\mu_N$ is appropriate to the
given voting system. Any of the empirical results of \cite{GelKB}
can be modeled by an appropriate choice of $\mu_N$.
\item It is only $\overline{\mu}_N$ that enters the formulae (\ref{res:geqsqrt}) and (\ref{res:leqsqrt}),
no other information about $\mu_N$ is relevant. The quantities $\overline{\mu}_N$ can be estimated using
Theorem \ref{th:weightmu}. In fact, more is true by the following result.
\end{enumerate}
\end{remarks}

\begin{theorem}\label{th:weakconv}
Let $P_N$ be the distribution of $\frac{1}{N}\,\sum_{i=1}^N\,X_i$ under the measure $\P_{\mu_N}$
then the sequence of measures
$P_N-\mu_N$ converges weakly to $0$.
\end{theorem}

Note that the distribution of $\frac{1}{N}\,\sum_{i=1}^N\,X_i$ is
the distribution of the voting results of the voter $i=1,\ldots,N$.
This is the quantity considered in \cite{GelKB}. Theorem
\ref{th:weakconv} tells us that the distribution of the voting
results for large number $N$ of voters is approximately equal to the
distribution $\mu_N$. In particular, for independent voting the
voting result is always extremely tight while for Straffin's example
any voting result has the same probability, i.e. it is equally
likely that a proposal gets 99\%  or 53\% of the votes.

\section{Voting models as spin systems\label{sec:spin}}
Spin systems are a central topic in statistical physics. They model
magnetic phenomena. The spin variables, usually denoted by $\sigma_i$, may take values in the set $\{-1, +1\}$
with $+1$ and $-1$ meaning `spin up' and `spin down' respectively. The spin variables model the elementary
magnets of the material (say the electrons or nuclei in a solid). The index $i$ runs over an index set $I$
which represents the set of elementary magnets.

The probability measure underlying the statistical structure is typically given by a `Gibbs measure'
defined through an energy functional $\mathcal{E}(\{\sigma_i\}_{i\in I})$. $\mathcal{E}$ gives the energy
of a given spin configuration $\{\sigma_i\}$. The system prefers configurations with low energy $\mathcal{E}$.
This is expressed in the Gibbs measure given by:

\begin{equation}\label{def:preGibbs}
q\,(\{\sigma_i\}_{i\in I}) = e^{-\beta\mathcal{E}(\{\sigma_i\}_{i\in
I})}.
\end{equation}

The parameter $\beta$ plays the role of an inverse temperature.
$q$ defines a (counting) measure on the space
$\Omega=\{-1,+1\}^I$. It has total mass:

\begin{equation}\label{def:Z}
   \mathcal{Z}\;=\;\sum_{\{\sigma_i\}\in\Omega}\,e^{-\beta\mathcal{E}(\{\sigma_i\}_{i\in
   I})}.
\end{equation}

Hence we obtain a \emph{probability} measure by setting:

\begin{equation}\label{def:Gibbs}
p\,(\{\sigma_i\}_{i\in I}) =
\mathcal{Z}^{-1}\;e^{-\beta\mathcal{E}(\{\sigma_i\}_{i\in I})}.
\end{equation}

Of course, we may interpret any spin system as a voting system with voting measure $p$,
as long as $\mathcal{E}(\{\sigma_i\})=\mathcal{E}(\{-\sigma_i\})$, and vice versa.

In particular, independent voting corresponds to the energy functional $\e(\{\sigma_i\})\equiv 1$.

Moreover, the \gf{common belief} model is given by an energy function:

\begin{equation}
\e(\{\sigma_i\})~=~-\;h\,\sum_i \,\sigma_i
\end{equation}

where $h$ is a random variable connected to the variable $Z$ defined in (\ref{def:distZ})
by:

\begin{equation}\label{eq:Zh}
\eh\, (1+Z) = \frac{e^h}{e^h+e^{-h}}.
\end{equation}

Note, that when $h$ runs from $-\infty$ to $\infty$ in (\ref{eq:Zh}) the value of $Z$ runs monotonously
from $-1$ to $+1$.

In term of statistical physics the above model is a system without spin-spin interaction in a random but
constant magnetic field. The inverse temperature $\beta$ which we encountered in equation (\ref{def:preGibbs})
is superfluous in this model as it can be absorbed in the magnetic field strength $h$.

\section{The voters' interaction model\label{sec:mfm}}
In the common belief model the voting behavior of each voter is
influenced by a preassigned, a priori given common belief variable
$Z$. The correlation between the voters results from the general
voting tendency described by the value of $Z$.

In this section we investigate a model with a \emph{direct}
interaction between the voters, namely a tendency of the voters to
vote in agreement with each other. In the view of statistical
physics this corresponds to the tendency of magnets to align. There
are various models in statistical physics to prescribe such a
situation. Presumably the best known one is the Ising model where
neighboring spins interact in the prescribed ways. The neighborhood
structure is most of the time given by a lattice (e.g. $\Z^d$). The
results on the system depend strongly on that neighborhood
structure, in the case of the lattice $\Z^d$ on the dimension $d$.

In the following we consider another, in fact easier model where no such assumption on the local
\gf{neighborhood} structure has to be made. We consider it an advantage of the model that very little
of the microscopic correlation structure of a specific voting system enters into the model.

The model we are going to consider is known in statistical mechanics as the \emph{Curie-Weiss model} or the
\emph{mean field model} (see e.g. \cite{Thompson}, \cite{BoltS} or \cite{Dorlas}). In
this model a given voter (spin) interacts with all the other voters (resp. spins) in a way which makes it
more likely for the voters (spins) to agree than to disagree. This is expressed through an energy function
$\e$ which is smaller if voters agree. Note that a \emph{small} energy for a given voting configuration
(relative to the other configurations) leads to a \emph{high} probability of that configuration
relative to the others through formula (\ref{def:Gibbs}).

The energy $\e$ for a given voting outcome $\{X_i\}_{i=1\ldots N}$ is given in the mean field model by:

\begin{equation}\label{def:meanf}
    \e\big(\{X_i\}\big)= ~-\,\frac{J}{N-1}\,\sum_{i,j\atop i\not=j}\;X_i
    X_j.
\end{equation}
Here $J$ is a non negative number called the coupling constant.
According to (\ref{def:meanf}) the energy contribution of a single
voter $X_i$ is expressed through the averaged voting result of all other
voters $\frac{1}{N-1}\,\sum_{j\not=i}X_j$. If $X_i$ agrees in sign
with this average the voter $i$ makes a negative contribution to the total
energy, otherwise $X_i$ will increase the total energy. The strength
of this negative or positive contribution is governed by the
coupling constant $J$. In other words: situations for which $X_i$
agrees with the other voters in average are more likely than others.
This can be seen from the formula for the probability of a given
voting outcome, namely:

\begin{equation}\label{def:pmf}
    p_J\,\big(\{X_i\}\big)\,=~\mathcal{Z}^{-1}\,e^{-\,\e(\{X_i\})}~=
    ~\mathcal{Z}^{-1}\,e^{\frac{J}{N-1}\,\sum_{i\not=j}\;X_i X_j}
\end{equation}
where we have set
\begin{equation}
    \mathcal{Z}~=~\sum_{\{X_i\}\in\,\{\pm1\}^N}\;e^{-\,\e(\{X_i\})}.
\end{equation}

As before the parameter $\beta$ is not needed, it can be absorbed in the coupling constant $J$.

Our goal is to compute the average:
\begin{equation}\label{optweimf}
    w_N~=~\E_{J,N}\big(\,|\sum_{i=1}^N\;X_i\,|\,\big).
\end{equation}

Here $\E_{J,N}$ denotes expectation with respect to the measure
defined in (\ref{def:pmf}). The quantity $w_N$ gives the optimal
weight in the council for a population of $N$ voters with a
correlation structure given by a mean-field model with coupling
constant $J$. We will see that the value of $w_N$ changes dramatically
when $J$ changes from a value below one to a value above one. This
has to do with the fact that the mean-field model undergoes a phase
transition at the point $J=1$ (see \cite{BoltS, Dorlas, Thompson}).

\begin{theorem}\label{th:mf}~\\[-3mm]
\begin{enumerate}
\item\label{th:mf1} If $J<1$ then
\begin{equation}\label{less1th:mf}
    w_N~=~\E_{J,N}\big(\,|\sum_{i=1}^N\;X_i\,|\,\big)~\sim~\frac{\sqrt{2}}{\sqrt{\pi}}\;\frac{1}{\,\sqrt{1-J}\,}\; \sqrt{N}\qquad\textnormal{as }
    N\to\infty.
\end{equation}
\item\label{th:mf2} If $J>1$ then
\begin{equation}\label{bigger1th:mf}
    w_N~=~\E_{J,N}\big(\,|\sum_{i=1}^N\;X_i\,|\,\big)~\sim~C(J)\; N\qquad\textnormal{as }
    N\to\infty.
\end{equation}
\end{enumerate}
\end{theorem}

\begin{remarks}~\\[-3mm]
\begin{enumerate}
\item By $x_N\sim y_N$ as $N\to\infty$ we mean that
$\lim_{n\to\infty}\frac{x_N}{y_N}=1$.

\item The constant $C(J)$ in (\ref{bigger1th:mf}) can be computed: If $J>1$ then $C(J)$ is the (unique)
positive solution $C$ of
\begin{equation}\label{mf:Ccritical} ~\tanh(J\,C)=C.
\end{equation}
Note that for $J\leq 1$ there is no positive solution of equation of
\ref{mf:Ccritical}.
\end{enumerate}
\end{remarks}

The proof of Theorem \ref{th:mf} will be given in section \ref{sec:pr2}.

\section{Conclusions\label{sec:conc}}
The above calculations show that one can reproduce the square-root
law as well as the results of \cite{GelKB} and other laws by
assuming particular correlation structures among the voters of a
certain country. To find the right model is a question of adjusting
the parameters of the models to empirical data of the country under
consideration. Moreover, the models allow us to investigate
questions about voting systems on a theoretical level. We believe that the models
described above can help to understand voting behavior in many situations.

To \emph{design} a nonhomogeneous voting system for a
\emph{constitution} in the light of our results is a question of
different nature. Even knowing the correlation structure of the
countries in question exactly would be of limited value to design a
constitution. Constitutions are meant for a long term period,
correlation structures of countries on the other hand are changing
even on the scale of a few years.

One might argue that modern societies have a tendency to decrease the correlation between their members.
In all modern states, at least in the West, the influence of churches, parties, and unions is constantly declining.

In addition to this it seems more important to protect small countries against a domination of the big ones than
the other way round. This motivates us to choose a square-root law in these long term cases.

\section{Appendix 1: Proofs for section \ref{sec:model}\label{sec:pr}}

We start with a short Lemma:

\begin{lemma}\label{lem:sym}
Suppose $X_1,...,X_N$ are $\{ -1,1 \} -$valued random variables with
the symmetry property (\ref{eq:symm}) then
\begin{eqnarray}\label{eq:mean}
\mathbb{E}(\sum_{i=1}^N X_i) = 0
\end{eqnarray}
and
\begin{eqnarray}\label{eq:abs}
\mathbb{E}(\sum_{i=1}^N X_i \ \chi(\sum_{i=1}^N X_i)) = \frac{1}{2}
\mathbb{E} (| \sum_{i=1}^N X_i |) \ .
\end{eqnarray}
\end{lemma}

\begin{remark}
As defined above $\chi (x)=1$ if $x>0$, $\chi (x)=-1$ if $x\leq 0$.
\end{remark}

\begin{proof}
(\ref{eq:symm}) implies
\begin{displaymath}
\mathbb{P}(X_i = 1) = \mathbb{P}(X_i = -1) = \frac{1}{2}
\end{displaymath}
hence $\mathbb{E}(X_i)=0$ and (\ref{eq:mean}) follows.

To prove (\ref{eq:abs}) we observe that due to (\ref{eq:symm})
\begin{eqnarray}\nonumber
\mathbb{E}(|\sum_1^N X_i|) &=& \mathbb{E}(\sum_{i=1}^N X_i \
\chi(\sum_{i=1}^N X_i)) - \mathbb{E}(\sum_{i=1}^N X_i \
\chi(-\sum_{i=1}^N X_i))
\\ \nonumber &=&2\mathbb{E}(\sum_{i=1}^N X_i \ \chi(\sum_{i=1}^N
X_i)) \ .
\end{eqnarray}
\end{proof}

We turn to the proof of Theorem \ref{th:srl}.

\begin{proof}(Theorem \ref{th:srl})
Let us abbreviate: $S_\nu :=\sum_{i=1}^{M_\nu}X_{\nu i}$.

Observe that the $S_\nu$ are independent by assumption and satisfy
$\mathbb{E}(S_\nu)=0$, moreover
\begin{eqnarray}\label{eq:S0}
\mathbb{E}(S_\nu \ \chi (S_\mu)) = 0 \textnormal{ if } \nu\neq\mu
\end{eqnarray}
and
\begin{eqnarray}\label{eq:sabs}
\mathbb{E}(S_\nu \ \chi (S_\nu)) = \frac{1}{2}\mathbb{E}(| S_\nu |)
\end{eqnarray}
by Lemma \ref{lem:sym}. To find the minimum of the function
\begin{displaymath}
\Delta(w_1 ,..., w_M ) = \mathbb{E}( (\sum_1^M S_\nu - \sum_1^M
w_\nu \ \chi(S_\nu))^2 )
\end{displaymath}
we look at the zeros of $\frac{\partial\Delta}{\partial w_\mu}$.
\begin{eqnarray}\nonumber
0=\frac{\partial\Delta}{\partial w_\mu} &=& -2\mathbb{E}\big((\sum_1^M
S_\nu - \sum_1^M w_\nu \ \chi(S_\nu))\chi(S_\mu)\big)
\\ \nonumber &=& -2\mathbb{E}(S_\mu\ \chi (S_\mu) - w_\mu \ \chi (S_\mu) \ \chi
(S_\mu)) \ .
\end{eqnarray}
So
\begin{displaymath}
w_\mu\; \mathbb{E}( (\chi (S_\mu))^2 ) = \mathbb{E}({S_\mu \ \chi
(S_\mu)}) = \eh\;\mathbb{E}(|S_\mu|) \ .
\end{displaymath}
Since $\chi (S_\mu)^2 = 1$ we obtain
\begin{displaymath}
w_\mu = \frac{1}{2}\;\mathbb{E}(|S_\mu|) \ .
\end{displaymath}
\end{proof}

We turn to the proof of Theorem \ref{th:indep}.

\begin{proof}
Let $X_1,...,X_N$ be $\{ -1,1 \} -$valued random variables with
\\$P(X_i = 1)= P(X_i = -1)=\frac{1}{2}$. Then
\begin{displaymath}
\mathbb{E}(| \sum_1^N X_i |) = ~\sqrt{N}~\; \mathbb{E}(|
\frac{1}{\sqrt{N}}\sum_1^N X_i |) \ .
\end{displaymath}
By the central limit theorem (see e.g. \cite{Lamperti})
$\frac{1}{\sqrt{N}}\sum_1^N X_i$ has asymptotically a normal
distribution with mean zero and variance 1, hence \\$\mathbb{E}(|
\frac{1}{\sqrt{N}}\sum_1^N X_i |)\rightarrow \frac{\sqrt{2}}{\sqrt{\pi}}$.
\end{proof}

\section{Appendix 2: Proofs for Section \ref{sec:cbm}\label{appii}}

\begin{proof}(Theorem \ref{th:corrpmu})
Since $\mathbb{E}_\mu(X_i)=0$,
\begin{eqnarray}
\langle X_i , X_j \rangle_\mu &=& \mathbb{E}_\mu (X_iX_j)
\\ \nonumber &=& \mathbb{P}_\mu
(X_i=X_j=1)+\mathbb{P}_\mu(X_i=X_j=-1)-2\mathbb{P}_\mu(X_i=1,
X_j=-1)
\\ \nonumber &=& \int d\mu (\zeta) \{ P_{\frac{1}{2}(1+\zeta)}(X_i=X_j=1)+ P_{\frac{1}{2}(1+\zeta)}(X_i=X_j=-1) \\
\nonumber &&\quad\quad\quad\quad  - 2P_{\frac{1}{2}(1+\zeta)}(X_i=1,
X_j=-1) \}
\\ \nonumber &=& \int d\mu (\zeta) \{ \frac{1}{4}(1+\zeta)^2 + \frac{1}{4} (1-\zeta)^2- \frac{1}{2}(1-\zeta^2) \}
\\ \nonumber &=& \int \zeta^2 d\mu (\zeta) \ .
\end{eqnarray}
\end{proof}

To prove Theorem \ref{th:weightmu}
we need the following Lemma:

\begin{lemma}\label{lem:centmean}
$\mathbb{E}_\mu(\frac{1}{N}|\sum (X_i-Z)|) \leq
\frac{1}{\sqrt{N}}$.
\end{lemma}

\begin{proof}
\begin{eqnarray}\nonumber
\mathbb{E}_\mu\big(\frac{1}{N}|\sum (X_i-Z)|\big) &=& \frac{1}{N}\,
\mathbb{E}_\mu\big(|\sum (X_i-Z)|\big)
\\ \nonumber &\leq& \frac{1}{N}\Big\{\;\mathbb{E}_\mu\Big(\big(\sum(X_i-Z  )\big)^2\Big)\Big\}^{1/2}
\\ \label{eq:var} &=& \frac{1}{N}\, \Big\{ \int d\mu (\zeta)\; E_{p_\zeta}\Big(\big( \sum_1^N (X_i-\zeta)
\big)^2\Big)\Big\}^{1/2}.
\end{eqnarray}
Given $Z=\zeta$ the random variables $X_i-\zeta$ have mean zero
and are independent with respect to the measure $P_{p_\zeta}$, thus
\begin{displaymath}
E_{p_\zeta}\Big( \big(\sum_1^N (X_i-\zeta)\big)^2 \Big) = N E_{p_\zeta}(X_i-\zeta)^2 =
N(1-\zeta^2)\leq N \ ,
\end{displaymath}
hence
\begin{displaymath}
(\ref{eq:var}) \leq \frac{1}{\sqrt{N}}(\int d\mu (\zeta)
(1-\zeta^2))^{1/2} \leq \frac{1}{\sqrt{N}} \ .
\end{displaymath}
\end{proof}

Using Lemma \ref{lem:centmean} we are in a position to prove
Theorem \ref{th:weightmu}:

\begin{proof}
(1) Suppose that:
\begin{equation}
    \overline{\mu}_N=\int |\zeta|\;d\mu_N(\zeta)\geq \,C\,\frac{1}{N^{1/2-\varepsilon}}
    \end{equation}
then we estimate:
\begin{eqnarray}\nonumber
\mathbb{E}_{\mu_N}(\frac{1}{N}| \sum_1^N X_i |) &=& \mathbb{E}_{\mu_N}(|
\frac{1}{N}\sum_1^N (X_i-Z) + Z |)
\\ \nonumber &\leq& \mathbb{E}_{\mu_N}(|Z|) + \mathbb{E}_{\mu_N}(| \frac{1}{N}\sum_1^N (X_i-Z) |)
\\ \label{estim1} &\leq& \overline{\mu}_N + \frac{1}{\sqrt{N}}
\end{eqnarray}
by Lemma \ref{lem:centmean}. Moreover
\begin{eqnarray}\nonumber
\mathbb{E}_{\mu_N}(\frac{1}{N}| \sum_1^N X_i |) &\geq&
\mathbb{E}_{\mu_N}(|Z|)-\mathbb{E}_{\mu_N}( | \frac{1}{N}\sum X_i-Z
| )
\\ \label{estim2} &\geq& \overline{\mu}_N - \frac{1}{\sqrt{N}} \
.
\end{eqnarray}

Hence
\begin{equation}
|\,\mathbb{E}_{\mu_N}(\frac{1}{N}| \sum_1^N X_i |)\;-\; \overline{\mu}_N\,|~\leq~\frac{1}{\sqrt{N}}
\end{equation}
which proves (\ref{res:geqsqrt}).

(2) To prove (\ref{res:leqsqrt}) we obtain by the same reasoning as above:
\begin{equation}
|\,\mathbb{E}_{\mu_N}(\frac{1}{N}| \sum_1^N X_i |)\;-\; \overline{\mu}_N\,|~\leq~\frac{1}{\sqrt{N}}
\end{equation}
\end{proof}

We end this section with the proof of Theorem \ref{th:weakconv}:

\begin{proof} We have to prove that for bounded continuous functions $f$:
\begin{equation}\label{eq:weakconv}
\int\,\Big(\,f(\frac{1}{N}\,\sum_{i=1}^N\,X_i)-f(Z)\,\Big)\;d\,\P_{\mu_N}~\to~0.
\end{equation}

The convergence (\ref{eq:weakconv}) is clear for continuously differentiable $f$ from Lemma \ref{lem:centmean}.
It follows for arbitrary bounded continuous $f$ by a density argument.
\end{proof}

\section{Appendix 3: Proofs for section \ref{sec:mfm}\label{sec:pr2}}
In this section we prove Theorem \ref{th:mf}.

\begin{proof} (Theorem \ref{th:mf} (\ref{th:mf1})\,)\\
We denote by $E_0^{(N)}$ the expectation of the coin tossing model
for $N$ independent symmetric $\{+1,-1\}$-valued random variables, i.e.:
\begin{equation}
E_0^{(N)}\big(F(X_1,\ldots,
X_N)\big)~=~\frac{1}{2^N}\,\sum_{\{x_i\}\in\,\{+1,-1\}^N}\,f(x_1,\ldots,x_N).
\end{equation}

\noindent We set:
\begin{equation}
\mathcal{Z}_{J N} = E_0^{(N)}\Big(\,e^{\frac{J}{2}\big(\frac{1}{\sqrt{N}}\sum_{i=1}^N X_i\big)^2}\Big)
\end{equation}
and:
\begin{equation}
\mathcal{X}_{J N} =
\;E_0^{(N)}\Big(\,|\frac{1}{\sqrt{N}}\sum_{i=1}^N X_i|\;
e^{\,\frac{J}{2}\big(\frac{1}{\sqrt{N}}\sum_{i=1}^N
X_i\big)^2}\Big).
\end{equation}

\noindent Then:
\begin{equation}\label{eq:absX}
    \E_{J N}(|\sum_{i=1}^N
    X_i|)~=~\sqrt{N}\;\;\frac{\mathcal{X}_{J,N}}{\mathcal{Z}_{J,N}}.
\end{equation}

\noindent Under the probability law $E_0^{(N)}$ the random variables $X_i$ are centered and independent,
thus the central limit theorem (see e.g. \cite{Lamperti}) tells us that $\frac{1}{\sqrt{N}}\sum_{i=1}^N X_i$
converges in distribution to a standard normal distribution. Consequently, for $J<1$ and $N\to\infty$ :

\begin{equation}
\mathcal{Z}_{J N} \to \frac{1}{\sqrt{2\pi}}\,\int_{-\infty}^{\infty}\;e^{\,-\,\frac{(1-J)\,x^2}{2}}\;dx~=
~\frac{1}{\,\sqrt{1-J\,}\,}
\end{equation}

and:
\begin{equation}
\mathcal{X}_{J N} \to
\frac{1}{\sqrt{2\pi}}\,\int_{-\infty}^{\infty}\;|x|\;e^{\,-\,\frac{(1-J)\,x^2}{2}}\;dx
~=~\frac{\sqrt{2}}{\sqrt{\pi}}\;\frac{1}{\,1-J\,}.
\end{equation}

Consequently:

\begin{equation}
    \E_{J N}(|\sum_{i=1}^N X_i|)~=~\sqrt{N}\;\;\frac{\mathcal{X}_{J,N}}{\mathcal{Z}_{J,N}}
    ~\sim~ \frac{\sqrt{2}}{\sqrt{\pi}}\;\frac{1}{\,\sqrt{1-J}\,}\;\sqrt{N}.
\end{equation}

\end{proof}

\begin{proof} (Theorem \ref{th:mf} (\ref{th:mf2})\,)\\
By Theorem 6.3 in \cite{BoltS} the distribution $\nu_N$ of $S_N=\frac{1}{N}\sum_{i=1}^N X_i$
converges weakly to the measure $\nu=\delta_{-C(J)}+\delta_{C(J)}$ where $C(J)$ was defined
in (\ref{mf:Ccritical}).

Hence,
\begin{eqnarray}
\E_J(|\sum_{i=1}^N X_i|) &=& N\; \E_J(|S_N|)~\\
&=&~N\;\int|\lambda|\;d\nu_N(\lambda)\\
&\approx& ~N\;\int|\lambda|\;d\nu(\lambda) \\
&=&~N\;C(J).
\end{eqnarray}
\end{proof}

 \end{document}